\documentclass[reqno,11pt]{amsart}

\input xy
\xyoption{all}
\usepackage{graphicx}
\usepackage{transparent}
\usepackage[mathscr]{eucal}
\usepackage{amsmath,amssymb,amsfonts,bbm,stmaryrd,hhline, amsthm}
\usepackage{amscd}
\usepackage{amsopn}
\usepackage{url}
\usepackage[breaklinks]{hyperref}
\hypersetup{colorlinks}
\usepackage[T1]{fontenc}

\usepackage{color} %\textcolor{red}{Test}

\definecolor{darkred}{rgb}{1,0,0} %can change the intensity in [0,1]
\definecolor{darkgreen}{rgb}{0,0.8,0}
\definecolor{darkblue}{rgb}{0,0,1}

\hypersetup{colorlinks, linkcolor=darkblue, filecolor=darkgreen,
urlcolor=darkred, citecolor=darkgreen}

\numberwithin{equation}{section}
\newtheorem {Theorem}{Theorem}
\numberwithin{Theorem}{section}

\newtheorem {Lemma}[Theorem]    {Lemma}

\newtheorem {Proposition}[Theorem]{Proposition}

\theoremstyle{definition}

\theoremstyle{remark}
\newtheorem{Remark}[Theorem]{Remark}
\newtheorem{Example}[Theorem]{Example}

\expandafter\chardef\csname pre amssym.def
at\endcsname=\the\catcode`\@ \catcode`\@=11
\def\undefine#1{\let#1\undefined}
\def\newsymbol#1#2#3#4#5{\let\next@\relax
 \ifnum#2=\@ne\let\next@\msafam@\else
 \ifnum#2=\tw@\let\next@\msbfam@\fi\fi
 \mathchardef#1="#3\next@#4#5}
\def\mathhexbox@#1#2#3{\relax
 \ifmmode\mathpalette{}{\m@th\mathchar"#1#2#3}%
 \else\leavevmode\hbox{$\m@th\mathchar"#1#2#3$}\fi}
\def\hexnumber@#1{\ifcase#1 0\or 1\or 2\or 3\or 4\or 5\or 6\or 7\or 8\or
 9\or A\or B\or C\or D\or E\or F\fi}

\font\teneufm=eufm10 \font\seveneufm=eufm7 \font\fiveeufm=eufm5
\newfam\eufmfam
\textfont\eufmfam=\teneufm \scriptfont\eufmfam=\seveneufm
\scriptscriptfont\eufmfam=\fiveeufm

\catcode`\@=\csname pre amssym.def at\endcsname

\newcommand{\MM}{\widetilde{M}}
\newcommand{\LL}{\widetilde{\mathcal{L}}}

\newcommand{\CA}{{\mathcal A}}

\newcommand{\CS}{{\mathcal S}}

\newcommand{\A}{{\mathcal A}}

\newcommand{\D}{{\Delta}}

\def    \reals      {{\mathbb R}}
\def    \Z      {{\mathbb Z}}
\def    \N      {{\mathbb N}}
\def    \Q      {{\mathbb Q}}
\def    \TT     {{\mathbb T}}
\def    \T      {{\mathbb T}}
\def    \CP     {{\mathbb C}{\mathbb P}}

\def    \12    {{\frac{1}{2}}}

\def    \p      {\partial}

\def    \im     {\operatorname{im}}

\def    \Sp     {\operatorname{Sp}}

\def    \MUCZ  {\operatorname{\mu_{\scriptscriptstyle{CZ}}}}

\def    \WFS  {\operatorname{\omega_{\scriptscriptstyle{FS}}}}

\begin{document}

%%%%%%%%%%%%%%%%%%%%%%%%%%%%%%
%   TEXT FORMATTING

%\setlength{\textheight}{7.5in} \setlength{\textwidth}{5in}
%\setlength{\smallskipamount}{6pt}
%\setlength{\medskipamount}{10pt}
%\setlength{\bigskipamount}{16pt}

%%%%%%%%%%%%%%%%%%%%%%%%%%

%%%%%%%%%%%           BEGINNING OF  TEXT

%%%%%%%%%%%%%%%%%%%%%%%%%%

\title[On hyperbolic points and periodic orbits of symplectomorphisms]{\small On hyperbolic points and periodic orbits of symplectomorphisms}

\author[Marta Bator\'eo]{Marta Bator\'eo}

\address{IMPA - Instituto Nacional de Matem\'atica Pura e Aplicada;
Estrada Dona Castorina 110\\
Rio de Janeiro 22460-320, Brazil}
\email{mbatoreo@impa.br}

\begin{abstract}
We prove the existence of infinitely many periodic orbits of symplectomorphisms isotopic to the identity if they admit at least one hyperbolic periodic orbit and satisfy some condition on the flux. Our result is proved for a certain class of closed symplectic manifolds and the main tool we use is a variation of Floer theory for symplectomorphisms, the Floer-Novikov theory. 

The proof relies on an important result on hyperbolic orbits, namely, that a \emph{Floer-Novikov} trajectory which converges to an iteration $k$ of the hyperbolic orbit and crosses its fixed neighborhood has energy bounded below by a strictly positive constant independent of $k$. The main theorem follows from this feature of hyperbolic orbits and certain properties of quantum homology on the class of symplectic manifolds we work with.
\end{abstract}

\maketitle
\let\thefootnote\relax\footnotetext{
The author was partially supported by the FCT/Portugal grant SFRH/BD/44125/2008, the NSF grant DMS-1007149 and IMPA.}

\tableofcontents
\section{Introduction and the Main Result}
\subsection{Introduction}
The main result of this paper is the existence of infinitely many periodic orbits of symplectomorphisms isotopic to the identity as long as they admit at least one hyperbolic periodic orbit and satisfy some constraints on the flux. This is established for a certain class of closed symplectic manifolds using a variation of Floer theory for symplectomorphisms, the Floer-Novikov theory (\cite{LO_fixedpts95, O_flux06}). 

Although here we are interested in symplectomorphisms which do not necessarily arise from Hamiltonians, the main result of this paper can be viewed in the context of the Conley conjecture (\cite{Co}) which claims the existence of infinitely many periodic orbits (of a Hamiltonian diffeomorphism). The conjecture holds true for many symplectic manifolds, namely, with ${c_1}|_{\pi_2(M)}=0$ and negative monotone manifolds; see \cite{CGG11, GG:action09, He12} and also \cite{FrHa03, Gi:conley, GG:conley12, Hi, LeC, SZ92}. The main difference between the Conley conjecture and our result is that in the Conley conjecture the existence of periodic orbits is unconditional whereas in our result the symplectomorphism is required to have one contractible (hyperbolic) periodic orbit. Without loss of generality, when a periodic orbit is contractible, we may assume it is a fixed point. Hence, for the sake of simplicity, from now on we consider the hyperbolic periodic orbit, say $\gamma$, to be constant; see beginning of Section~\ref{proof:maintheo} for more details.

Due to this assumption on the existence of a fixed point of a specific type, our result fits more accurately under what G{\"u}rel describes in \cite{Gu_noncon, Gu_linear} as the generalized HZ-conjecture; see also \cite{GG:hyp12}. This variant of the Conley conjecture claims that a Hamiltonian diffeomorphism with ``more than necessary" fixed points has infinitely many periodic points. Here ``more than necessary" is understood as a lower bound provided by some form of the Arnold conjecture, e.g., the expected bound for $\CP^n$ is $n+1$. The HZ-conjecture was originally stated (as far as we know) in this form by Hofer and Zehnder in \cite[p.263]{HZ11} and was motivated by the results of Gambaudo and Le Calvez in \cite{Gam_LeCa} and Franks in \cite{Fr88} (see also \cite{Fr92, Fr96}) where they prove that
\vspace*{0.1cm}
\begin{center}
\emph{an area preserving diffeomorphism of $S^2$
with at least three fixed points\\
 has infinitely many periodic points};
\end{center} 
\vspace*{0.1cm}
see also \cite{BH11, CKRTZ12, Ke12} for symplectic topological proofs.
A broader variant of the conjecture suggests the presence of infinitely many periodic points provided the existence of a fixed point which is \emph{unnecessary} from a homological or geometrical perspective. In fact, our main theorem (Theorem~\ref{maintheo}) asserts that, for a certain class of symplectic manifolds, a symplectomorphism (symplectically isotopic to the identity) with a hyperbolic fixed point must admit infinitely many periodic points (as long as it satisfies some condition on its flux). The theorem is a symplectic analogue of a result proved in \cite{GG:hyp12} for Hamiltonian diffeomorphisms. In dimension greater than two, the conjecture is also supported by the result presented in \cite{Gu_noncon} for non-contractible orbits and by its ``local version'' considered in \cite{Gu_linear}.

%Among the manifolds meeting the requirements of Theorem~\ref{maintheo} are complex projective spaces $\CP^n$, complex Grassmannians $\text{Gr} (2,N)$, $\text{Gr}(3,6)$ and $\text{Gr}(3,7)$, the products $\CP^n\times P^{2m}$  (with $m\leq n$) and $\text{Gr}(2,N)\times P^{2m}$ (with $m\leq 2$) where $P$ is symplectically aspherical and monotone products $\CP^n\times\CP^n$. Recall that we are interested in symplectomorphisms which are not Hamiltonian diffeomorphisms (see \cite{GG:hyp12} for the Hamiltonian case). In some instances, e.g. $\CP^n$, a symplectomorphism symplectically isotopic to the identity is automatically a Hamiltonian diffeomorphism. Some examples of manifolds which meet the requirements of the main theorem and admit symplectomorphisms (which are not Hamiltonian diffeomorphisms) with periodic points are the products of complex projective spaces with tori $\CP^n\times\TT^{2m}$ (with $m\leq n$) and the products of complex Grassmannians with tori $\text{Gr}(2,N)\times\TT^{2m}$ (with $m\leq 2$).

Some examples of manifolds which meet the requirements of the main theorem and admit symplectomorphisms (which are not Hamiltonian diffeomorphisms) with periodic points are the products of complex projective spaces with symplectically aspherical manifolds $\CP^n\times P^{2m}$ (with $m\leq n$) and the products of complex Grassmannians with symplectically aspherical manifolds $\text{Gr}(2,N)\times P^{2m}$ (with $m\leq 2$); see Examples~\ref{example:firstchern}~and~\ref{example:quantumprod} for further details.

\subsection{Existence of Infinitely Many Periodic Orbits (Main Theorem)} Consider a symplectomorphism $\phi$ in the identity component of the group of symplectomorphisms of $(M,\omega)$. 
%The flux homomorphism associates with $\phi$, or rather with an element of the universal covering of the group, a cohomology class $[\theta]$ in $H^1(M,\reals)$ (see definition in Section~\ref{section:sympmrpfisms}).
%We say that $\phi$ has \emph{rational flux} 
The flux homomorphism $\widetilde{\text{Flux}}$ associates with an element $[\phi_t]$ of the universal covering of $Symp_0(M,\omega)$ a cohomology class $[\theta]\in H^1(M;\reals)$. We say that the symplectic isotopy $\phi_t$ has \emph{rational flux} if the group formed by the integrals of $\theta$ over the loops in $M$ is
discrete, that is,
$$
\left< [\theta], \pi_1(M) \right>=h_1 \Z
$$
for some $h_1\in \reals.$

Then, we have the following result on the periodic points of such symplectomorphisms.

\begin{Theorem}[Main Theorem]\label{maintheo}
\

Let $M^{2n}$ be closed and strictly monotone (i.e. $M$ is monotone and $c_1|_{\pi_2(M)} \not=0$ and $[\omega]|_{\pi_2(M)}\not=0$). Assume that
\begin{itemize}
\item $c_1^{\min}\geq n/2 +1$ (where $c_1^{\min}$ is the minimal Chern number) and
\item[] \vspace*{-0.5cm}\begin{eqnarray}\label{l:1}
\hspace*{-2.55cm}\bullet \hspace*{0.15cm} \beta*\alpha=\mathfrak{q} [M]\footnotemark \quad\mbox{ in } HQ_*(M)=H_*(M)\otimes \Lambda
\end{eqnarray}
for some ordinary homology classes $\alpha,\;\beta\in H_*(M)$\\ 
(with $\deg(\alpha),\;\deg(\beta)<2n$).
\end{itemize}
Then any symplectomorphism $\phi$ in $\mbox{Symp}_0(M,\omega)$ with
\begin{itemize}
\item a contractible hyperbolic periodic orbit $\gamma$ and
\item $h_1/h_0\in \Q$, where $h_0$ is the rationality constant of $(M, \omega)$ and $h_1$ the rationality constant of the flux of a symplectic isotopy $\phi_t$ such that $\phi_0=id$ and $\phi_1=\phi$ 
%rational flux where $h_1=(p/r)h_0$ for some isotopy from the identity to $\phi$ (where $h_0$ is the rationality constant of $M$)
\end{itemize}
has infinitely many periodic orbits.

Here $\mathfrak{q}$ is the element of the Novikov ring, with degree $-2c_1^{\min}$, defined as in Section~\ref{subsubsection:novikov}.
\end{Theorem}

The assumption on the existence of a hyperbolic periodic orbit $\gamma$ is extremely important. A significant feature of these orbits is the fact that the energy needed for a (\emph{Floer-Novikov}) trajectory to converge to a $k$-th iteration of $\gamma$ and crossing its fixed neighborhood is bounded below by a strictly positive constant independent of $k$ (see Section~\ref{section:BCEthm}). The main tool used to prove our result is filtered Floer-Novikov homology (see Section~\ref{section:fnhomol}).

The proof of the theorem goes by contradiction and, if a symplectomorphism admits only finitely many periodic points, there exists an iteration $k$ for which the \emph{action} of all the $k$-periodic orbits are in a small neighborhood of $h_0\Z$. The assumption that the flux and the action grow together, more precisely $h_1/h_0\in \Q$, plays an essential role in this part of the argument. In this case, the action spectrum is not dense and it can be shown that the action values of the $k$-periodic orbits are, indeed, close to $h_0\Z$. Then using the feature of the hyperbolic orbit, the fact that quantum homology acts on the (filtered) Floer-Novikov homology and condition (\ref{l:1}), we obtain a $k$-periodic orbit with action outside the small neighborhood of $h_0\Z$. 
\begin{Remark} The hyperbolicity condition is required so that the orbit has the important feature mentioned above and which is also described in Section~\ref{section:BCEthm}. Hence, more \emph{generally}, a symplectomorphism with a periodic orbit with the  property described in Theorem~\ref{thm:ballcrossing} and satisfying the requirement on its flux admits infinitely many periodic points.\qed
\end{Remark}

The following proposition (proved in Section~\ref{section:prop}) leads to examples of symplectomorphisms which meet the requirements of the main theorem. 

In order to state the proposition, we first introduce the notion of \emph{Hamiltonian deformation}. We say that a symplectic isotopy $\phi_t$ of $(M,\omega)$ is obtained by a Hamiltonian deformation of another symplectic isotopy $\psi_t$ when $\phi_t\circ\psi_t^{-1}$ is a Hamiltonian diffeomorphism for all $t$. Notice that the flux is invariant by Hamiltonian perturbations.

\begin{Proposition}\label{prop:example}
Given a symplectic isotopy $\psi_t$ of $(M,\omega)$ and a contractible loop $\gamma$ in $M$, there exists an isotopy $\phi_t$ obtained by a Hamiltonian deformation of $\psi_t$ such that $\gamma$ is a hyperbolic one-periodic orbit of $\phi_1$.
\end{Proposition}

\begin{Example} Consider $M=\CP^n\times\TT^{2}$ with the standard symplectic form. Recall that $c_1(\TT^{2})=0$ and, under the normalization $\WFS[g_0]=n+1$ (where $\WFS$ is the Fubini-Studi form on $\CP^n$ and $g_0$ is the generator of $H_2(\CP^n,\reals)$), $\CP^n$ is a monotone symplectic manifold. Denote the rationality constant of $M$ by $h_0$. Recall that in $\T^2$ the flux of the isotopy
$$
\overline{\psi}_t(x,y)=(x+t\theta,y)\quad\quad\text{with}\quad \theta\not=1 
$$
is determined by the shift $\theta$, namely, $\widetilde{\text{Flux}}([\overline{\psi}_t])(y)=\theta$ where we identify $H^1(M,\reals)=\text{Hom}(H_1(M,\Z),\reals)$ and $y$ is the loop $\reals/\Z \rightarrow \T^2$ defined by $s\mapsto (0,s)$. Hence, the rationality constant $h_1$ of the flux (defined by (\ref{def:h_1})) is also determined by $\theta$. Now, consider the isotopy $\psi_t$ in $M$ defined by the identity on $\CP^n$ and $\overline{\psi}_t$ on $\T^2$. For the fixed $h_0$ in $(M,\omega)$, consider $\theta$ so that $h_1/h_0\in \Q$.
The symplectomorphism $\psi_1$ has no fixed points. According to the previous proposition, there exists an isotopy $\phi_t$ obtained by Hamiltonian deformation of $\psi_t$ such that $\phi_1$ has a hyperbolic fixed point.\qed
\end{Example}

\subsection{Acknowledgments} 
The author is grateful to Viktor Ginzburg for posing this problem and valuable discussions and to Leonardo Macarini for useful comments and remarks.

\section{Preliminaries}\label{section:prelim}
\subsection{Symplectic manifolds}
\label{section:sympmnfd}
Let $(M^{2n},\omega)$ be a \emph{closed} (i.e. compact and with no boundary) rational symplectic manifold and consider an almost complex
structure $J$ on $TM$ compatible with $\omega,$ i.e., such that $\left< \xi, \eta \right> := \omega(\xi, J\eta)$ is a Riemannian metric on $M.$
Recall that $(M, \omega)$ is said to be \emph{(spherically) rational} if the
group formed by the integrals of $\omega$ over the spheres in $M$ is
discrete, that is,
$$
\left< [\omega], \pi_2(M) \right>=h_0 \Z
$$
where $h_0\geq 0$. The constant $h_0$ is called the \emph{rationality constant}.

Since the space of almost complex structures compatible with $\omega$ is connected, the first Chern class $c_1\in H^2(M,\Z)$ is uniquely determined by $\omega$.
The \emph{minimal Chern number} of a symplectic manifold $(M,
\omega)$ is the positive integer $c_1^{\min}$ which generates the discrete group formed by the integrals
of $c_1$ over the spheres in $M$, i.e.,
$$
\left<c_1, \pi_2(M)\right> = c_1^{\min} \Z
$$
where $c_1^{\min}\in \Z^+$.

A symplectic manifold $(M,\omega)$ is called \emph{monotone} if the cohomology classes $c_1$ and
$[\omega]$ satisfy the condition
$$
[\omega]|_{\pi_2(M)} = \lambda \; {c_1}|_{\pi_2(M)}
$$
for some non-negative constant $\lambda\in\reals.$ If $\lambda\not=0$ and $c_1|_{\pi_2(M)}\not=0$, we say that $M$ is \emph{strictly monotone}.\\

\begin{Example}\label{example:firstchern}
Consider $M=\CP^n\times P^{2m}$ where $P^{2m}$ is symplectically aspherical (with $m\leq n$). We have  $c_1^{\min}(\CP^n)=n+1$. Hence, the minimal Chern number of 
$M$ is $c_1^{\min}(M)=n+1$ (see, e.g., \cite{MS12}  for more details). In Theorem~\ref{maintheo}, the assumption on the minimal Chern number
$$
c_1^{\min}(M)\geq \dim(M)/4 +1,
$$
in this case, is
$$
n+1\geq \frac{n+m}{2} +1,
$$
which is equivalent to $m\leq n$. \qed 
\end{Example}

\subsection{Symplectomorphisms}
\label{section:sympmrpfisms}
We denote by $\mbox{Symp} (M,\omega)$ the symplectomorphism group of $(M,\omega)$ and by $\mbox{Symp}_0 (M,\omega)$ the identity component in $\mbox{Symp} (M,\omega)$.

Let $\phi\in \mbox{Symp}_0 (M,\omega)$ and consider $\phi_t$ a symplectic path connecting the identity to $\phi$, i.e., $\phi_0=id$ and $\phi_1=\phi$. A vector field $X_t$ is defined by:
$$
\frac{d}{dt} \phi_t = X_t \circ \phi_t.
$$

Recall that the \emph{flux homomorphism} is defined on the universal covering of $\mbox{Symp}_0 (M,\omega)$ as follows:
\begin{eqnarray*}
\widetilde{\mbox{Flux}}: \widetilde{\mbox{Symp}_0} (M,\omega) &\rightarrow& H^1(M,\reals)\\
\text{[}\phi_t \text{]} & \mapsto & \left[ \displaystyle\int_0^1 \omega(X_t, \cdot ) dt \right].
\end{eqnarray*}

Let $\theta$ be a closed one-form which represents the cohomology class $\widetilde{\mbox{Flux}}([\phi_t])$.
Throughout this paper, we assume that the group formed by the integrals of $\theta$ over the loops in $M$ is
discrete, that is,
\begin{eqnarray}\label{def:h_1}
\left< [\theta], \pi_1(M) \right>=h_1 \Z
\end{eqnarray}
for some $h_1\in \reals$.
In this case, we say that $\phi$ has \emph{rational flux}.\\

Recall that given a time dependent Hamiltonian $H\colon S^1 \times M \rightarrow \reals $ (where $S^1=\reals/\Z$) the \emph{Hamiltonian vector field}, $X_H$, is defined by
$$
\omega(X_H, \cdot)=-dH.
$$
The time-dependent Hamiltonian flow $\varphi^t_H$ is given by
$$
\frac{d}{dt} \varphi^t_H = X_H \circ \varphi^t_H,
$$
and the \emph{time-one map} $\varphi^1_H$ is denoted by $\varphi_H$. A symplectomorphism $\phi$ is called a \emph{Hamiltonian diffeomorphism} if it is the time-one map of some Hamiltonian, i.e., $\phi=\varphi_H$. 
A Hamiltonian of period one $H\colon S^1 \times M \rightarrow \reals $ regarded as a Hamiltonian of period $k$ is denoted by $H^{\natural k}$.

\subsubsection{Tailed-capped loops}\label{section:tailedcapped} Let $\mathcal{L}M$ be the space of contractible loops in $M$ and $\Omega M$ be the space of based contractible loops in $M$. The map $ ev\colon \mathcal{L}M \rightarrow M $ defined by $x \mapsto x(0)$ is a fibration (see e.g.~\cite[pp~83]{Hu} for the details). It induces a long exact sequence on the homotopy groups and part of it is given by
$$
\pi_1(\Omega M) \rightarrow \pi_1(\mathcal{L}M) \rightarrow \pi_1(M).
$$
Since this fibration admits a section consisting of constant loops,

$$\pi_1(\mathcal{L}M) \cong \pi_1(\Omega M) \oplus \pi_1(M).$$ 
With the identification $\pi_1(\Omega M)\equiv \pi_2(M)$ (see e.g.~\cite[pp.5--7]{Adams78} for the details), we have
\begin{eqnarray}\label{eqn:loopsp}
\pi_1(\mathcal{L}M) \cong \pi_2(M) \oplus \pi_1(M).
\end{eqnarray}
Consider the homomorphisms 
$$
\overline{I_{\theta}},\;\overline{I_{c_1}},\;\overline{I_{\omega}}\colon \pi_1(\mathcal{L}M)\rightarrow \reals
$$
induced, respectively, by
$$
\begin{array}{cccc}\label{eqn:I_omega}
     I_{\theta} \colon \pi_1(M) \rightarrow \reals,&     \gamma  &  \mapsto  & \int_{\gamma} \theta,\\
     &&&\\
     I_{c_1}    \colon \pi_2(M) \rightarrow \reals,&     A       &  \mapsto  & -2\int_{A} c_1,      \\
     &&&\\ 
     I_{\omega} \colon \pi_2(M) \rightarrow \reals,&     A       &  \mapsto  & -\int_{A} \omega 
\end{array}
$$
and the covering space $\LL\MM$ of $\mathcal{L}M$ associated with the induced homomorphisms, i.e. the covering whose fundamental group is the subgroup $\ker\overline{I_{\theta}}\cap \ker\overline{I_{c_1}} \cap \ker\overline{I_{\omega}}$ of $\pi_1(\mathcal{L}M)$ (cf.~\cite[Remark~3.2]{O_flux06}). The deck transformation group of $\LL\MM \rightarrow \mathcal{L}M$ is 
\begin{eqnarray}\label{eqn:deckcovering}
\frac{\pi_1(\mathcal{L}M)}{\ker\overline{I_{\theta}}\cap \ker\overline{I_{c_1}} \cap \ker\overline{I_{\omega}}}.
\end{eqnarray}

Alternatively (cf.~\cite[pp~158]{LO_fixedpts95}), the covering $\LL\MM \rightarrow \mathcal{L}M$ can be constructed so that the following diagram commutes:
\begin{equation}\label{diagram:covering_space}
\xymatrix{
{\LL\MM} \ar[r]^{\widetilde{j}}\ar[d]^{\widetilde{\Pi}}& {\mathcal{L}\MM} \ar[r]^{\widetilde{ev}} \ar[d]^{\Pi} & {\MM} \ar[d]^{\pi} \\
{\LL M}             \ar[r]^{j}                & {\mathcal{L}M}             \ar[r]^{ev}             & {M}.
}
\end{equation}

Here $\MM$ is the covering space of $M$ associated with the subgroup $\ker I_{\theta}$ of $\pi_1(M)$. The deck transformation group of the covering $\pi: \MM \rightarrow M$ is isomorphic to
$$
\frac{\pi_1(M)}{\ker I_{\theta}}.
$$
The map 
$ev:\mathcal{L}M\rightarrow M$ is defined as above and the map $j\colon \LL M \rightarrow \mathcal{L}M$ is the covering associated with the homomorphisms $I_{c_1}$ and $I_{\omega}$. The deck transformation group of the covering $j$ is isomorphic to the quotient group
$$
\frac{\pi_2(M)}{\ker I_{c_1}\cap \ker I_{\omega}}.
$$

Hence, the deck transformation group of the covering $\LL\MM \rightarrow \mathcal{L}M$ is the direct sum
\begin{eqnarray}\label{eqn:decktransformation}
\frac{\pi_1(M)}{ \ker I_{\theta}} \oplus \frac{\pi_2(M)}{ \ker I_{c_1}\cap \ker I_{\omega}}
\end{eqnarray}
and~(\ref{eqn:deckcovering}) is isomorphic to~(\ref{eqn:decktransformation}).

An element of the covering space $\LL\MM$ is represented by an equivalence class of pairs $(\widetilde{x},\widetilde{v})$ where
\begin{itemize}
\item[i)] $\widetilde{x}$ is a loop in $\widetilde{M}$,
\item[ii)] $\widetilde{v}$ is a disk in $\widetilde{M}$ bounding $\widetilde{x}$ and
\item[iii)] $(\widetilde{x},\widetilde{v})$ is equivalent to $(\widetilde{y},\widetilde{w})$ if $\widetilde{x}=\widetilde{y}$ and
$$I_{c_1} (v\#(-w))=0=I_{\omega} (v\#(-w))$$ where $v=\pi\circ\widetilde{v}$ and $w=\pi\circ\widetilde{w}$.
\end{itemize}

\begin{figure}[htb]
  \centering
  \def\svgwidth{200pt}
  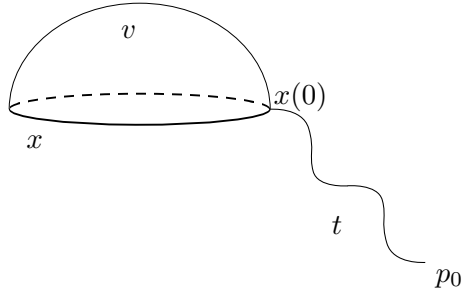
  \caption{\emph{Tailed}-capped loop}\label{figure:tc_loop}
\end{figure}

An element of $\LL\MM$ can be viewed as a capped loop with a \emph{tail} attached to it in $M$ (see Figure \ref{figure:tc_loop}): consider an element $[(\widetilde{x}, \widetilde{v})]\in \widetilde{\mathcal{L}}\widetilde{M}$ and a capped loop $(x,v)$ in $M$ (i.e. $x\colon S^1 \rightarrow M$ and 
$v\colon D^2 \rightarrow M$ satisfying $v|_{\partial D^2}=x$) such that $\pi\circ\widetilde{x}=x$ and $\pi\circ\widetilde{v}=v$. Fix a point $p_0\in M$ and consider a path, $t$, connecting $p_0$ to $x(0)$. We say that two objects $\widehat{\widehat{x}}:=(x,v,t)$ and $\widehat{\widehat{y}}:=(y,w,t^{\prime})$ are equivalent if
\begin{itemize}
\item[i)] $x=y$,
\item[ii)] $I_{c_1} (v\#(-w))=0=I_{\omega} (v\#(-w))$ and
\item[iii)] $I_{\theta}(t\#t^{\prime})=0$ where $t\#t^{\prime}$ is the concatenation of the paths $t$ and $t^\prime$.
\end{itemize}
The equivalence class $[{\widehat{\widehat{x}}}]$ (in $M$) corresponds to the equivalence class $[(\widetilde{x},\widetilde{v})]$ (in $\widetilde{M}$).

\subsubsection{Periodic orbits and the action functional}\label{section:po_af}
Let $\phi\in \mbox{Symp}_0(M)$, $\phi_t$ be a symplectic path connecting $\phi_0=id$ to $\phi_1=\phi$ and $X_t$ be the vector field associated with $\phi_t$ (as in Section \ref{section:sympmrpfisms}).

L\^e and Ono proved in \cite[Lemma~2.1]{LO_fixedpts95} that we can deform $\{\phi_t\}$ through symplectic isotopies (keeping the end points fixed) so that the cohomology classes $[\omega(X^{\prime}_t,\cdot)]$ (for all $t\in[0,1]$) and $\widetilde{\mbox{Flux}}([\phi_t])=[\theta]$ are the same (where $X^{\prime}_t$ is the vector field associated with the deformed symplectic isotopies $\phi^{\prime}_t$).
\begin{Lemma}[Deformation Lemma]\label{lemma:deformation}
%For $\{\phi_t\}\in \widetilde{\mbox{Symp}_0}(M,\omega)$, there exist
%\begin{itemize}
%\item[i)] a smooth path $\phi^{\prime}_t$ in $\mbox{Symp}_0(M,\omega)$ with $\phi^{\prime}_0=id$, $\phi^{\prime}_1=\phi_1$ and $\phi^{\prime}_{t+1}=\phi^{\prime}_t\circ\phi^{\prime}_1$ and
%\item[ii)] a Hamiltonian $H_t:M \rightarrow M$ with $H_{t+1}=H_t$
%\end{itemize}
%such that
%$$
%-\omega(X^{\prime}_t,\cdot)= \theta + dH_t=\colon\theta_t.
%$$
Each element in $\widetilde{Symp}_0(M,\omega)$ admits a representative symplectic isotopy generated by a smooth path of closed one-forms $\theta_t$ on $M$ whose cohomology class is identically equal to the flux.
\end{Lemma}

The fixed points of $\phi=\phi_1$ are in one-to-one correspondence with one-periodic solutions of the differential equation
\begin{eqnarray}\label{eqn:de_theta}
\dot{x}(t)=X_{\theta_t}(t,x(t))
\end{eqnarray}
where $X_{\theta_t}$ is defined by $-\omega(X_{\theta_t},\cdot)=\theta_t$.

The set of one-periodic solutions of (\ref{eqn:de_theta}) is denoted by $\mathcal{P}(\theta_t)$ and coincides with the zero set of the closed one-form defined on the loop space of $M$, $\mathcal{L}M$, by
\begin{eqnarray}\label{eqn:mv-af}
\alpha_{\{\phi_t\}}(x,\xi)= \displaystyle\int_0^1 \omega(\dot{x}, \xi) + \theta_t(x(t))(\xi) dt
\end{eqnarray}
where $x\in \mathcal{L}M$ and $\xi\in T_x\mathcal{L}M$ (i.e. $\xi$ is a tangent vector field along the loop $x$).

By the Deformation Lemma \ref{lemma:deformation}, there exists a periodic Hamiltonian
\begin{eqnarray}\label{eqn:H_tilde}
\widetilde{H}\colon S^1 \times \widetilde{M} \rightarrow \reals\quad \mbox{such that}\quad d\widetilde{H}_t=\pi^*\theta_t \quad (t\in S^1)
\end{eqnarray}
where $\theta_t:=-\omega(X_t,\cdot)$. The time-dependent Hamiltonian flow on $\widetilde{M}$ generated by $\widetilde{H_t}$ is the pull back of the original symplectic flow on $M$. In particular, the set of contractible periodic solutions of the Hamiltonian system $\mathcal{P}(\widetilde{H})$ is the set $\Pi^{-1}(\mathcal{P}(\theta_t))$ and $\mathcal{\widetilde{P}}(\widetilde{H}):=\widetilde{j}^{-1}(\mathcal{P}(\widetilde{H}))$ is the critical set of the functional:
\begin{eqnarray}\label{eqn:action_flnov}
\mathcal{A}_{\widetilde{H}}([\widetilde{x},\widetilde{v}])=-\int_D v^*\omega +\displaystyle\int_0^1 \widetilde{H}_t(\widetilde{x}(t)) dt 
\end{eqnarray}
where $\pi\circ\widetilde{v}=v$ (recall that $\widetilde{j}$ is given by (\ref{diagram:covering_space})). 
The action functional is homogeneous with respect to iterations in the following sense
$$
\mathcal{A}_{\widetilde{H}^{\natural k}}([\widetilde{x},\widetilde{v}]^{k})=k\mathcal{A}_{\widetilde{H}}([\widetilde{x},\widetilde{v}])
$$
where $[\widetilde{x},\widetilde{v}]^{k}$ is the $k$-th iteration of $[\widetilde{x},\widetilde{v}]$ and depends on the equivalence class of the capping $\widetilde{v}$ of the loop $\widetilde{x}$:
$$
\CA_{\widetilde{H}}([\widetilde{x},\widetilde{v}]\#A)=\CA_{\widetilde{H}}([\widetilde{x},\widetilde{v}]) + I_{\omega} (A)
$$
where $A\in \pi_2(M)$.

\subsubsection{Mean index and the augmented action}
Let $x$ be a periodic orbit of $\phi\in \mbox{Symp}_0(M,\omega)$, i.e., a periodic solution of (\ref{eqn:de_theta}). If the eigenvalues of the map  
$$
d\phi_{x(0)}:T_{x(0)}M\rightarrow T_{x(0)}M
$$
are not equal to one, the orbit $x$ is called \emph{non-degenerate}. If, in addition, none the of eigenvalues of the linearized return map $d\phi_{x(0)}$ has absolute value equal to one, we say $x$ is \emph{hyperbolic}. The symplectomorphism $\phi$ (or $H$ when $\phi=\varphi_H$) is said to be non-degenerate if all its one-periodic orbits are non-degenerate.

Let $(x,v)$ be a capped periodic orbit of $\phi$. Using a trivialization of $x^*TM$ arising from the capping $v$,
the linearized flow along $x$
$$
d\phi_t\colon T_{x(0)}M \rightarrow T_{x(t)}M
$$
can be viewed as a symplectic path $\Phi\colon[0,1] \rightarrow \Sp(2n).$ The mean index of $(x,v)$ is defined by $\D_{\phi}(x,v)
:=\D_{\phi}(\Phi)$; see \cite{SZ92}.

Recall that, in the setting of Section \ref{section:po_af}, the time-dependent Hamiltonian flow on $\widetilde{M}$ generated by $\widetilde{H}$ is the pull back of the original symplectic flow on $M$, hence a periodic orbit $\widetilde{x}\in \mathcal{\widetilde{P}}(\widetilde{H})$ is non-degenerate if and only if $\pi\circ\widetilde{x}$ is non-degenerate as a periodic orbit of $\phi$ and it is hyperbolic if and only if $\pi\circ\widetilde{x}$ is hyperbolic as a periodic orbit of $\phi$.  Moreover,
$$
\D_{\widetilde{H}}([\widetilde{x},\widetilde{v}])=\D_{\phi}((x,v))
$$
and it has the following properties:
\begin{eqnarray}\label{eqn:mean_homo}
\D_{\widetilde{H}^{\natural k}}([\widetilde{x},\widetilde{v}]^{k})=k\D_{\widetilde{H}}([\widetilde{x},\widetilde{v}]),
\end{eqnarray}
and
\begin{eqnarray}\label{eqn:mean_recap}
\D_{\widetilde{H}}([\widetilde{x},\widetilde{v}]\#A)=\D_{\widetilde{H}}([\widetilde{x},\widetilde{v}])+ I_{c_1}(A)
\end{eqnarray}
where $A\in \pi_2(M)$.
The \emph{augmented action} is defined by
\begin{eqnarray}\label{eqn:augmentedaction}
\underline{\mathcal{A}}_{\widetilde{H}}([\widetilde{x},\widetilde{v}]):= \mathcal{A}_{\widetilde{H}}([\widetilde{x},\widetilde{v}]) - \frac{\lambda}{2} \D_{\widetilde{H}}([\widetilde{x},\widetilde{v}]).
\end{eqnarray}
Notice that the augmented action is independent of the capping $\widetilde{v}$ and it is homogeneous with respect to iterations in the following sense
$$
\underline{\mathcal{A}}_{\widetilde{H}^{\natural k}}([\widetilde{x},\widetilde{v}]^{k})=k\underline{\mathcal{A}}_{\widetilde{H}}([\widetilde{x},\widetilde{v}]).
$$

\subsection{Floer-Novikov Homology}\label{section:fnhomol} In this section, we recall the construction of the Floer homology for symplectomorphisms following \cite{LO_fixedpts95} (references therein and \cite{O_flux06}).

Let $\phi$ be a symplectomorphism isotopic to the identity defined on a strictly monotone manifold $M$ and $\{\phi_t\}$ a path connecting $\phi$ to the identity. Take an almost complex structure $J$ on $M$ and fix an almost complex structure $\widetilde{J}$ on $\widetilde{M}$ corresponding to $J$. Consider the Hamiltonian $\widetilde{H}$ associated with $\phi_t$ as in (\ref{eqn:H_tilde}) and recall that we denote by ${\mathcal{P}}(\widetilde{H})$ the set of contractible periodic orbits of the Hamiltonian system associated with $\widetilde{H}$ and $\widetilde{\mathcal{P}}(\widetilde{H}):=\widetilde{j}^{-1}({\mathcal{P}}(\widetilde{H}))=$ Crit $(\mathcal{A}_{\widetilde{H}})$ (see Section~\ref{section:po_af}). The maps $\widetilde{u}\colon \reals \times S^1 \rightarrow \widetilde{M}$ which satisfy the equation
\begin{eqnarray}\label{eqn:connectorbit_eq_tilde}
\partial_s \widetilde{u} + \widetilde{J}(\widetilde{u})(\partial_t \widetilde{u} - X_{\widetilde{H}}(\widetilde{u}))=0
\end{eqnarray}
with boundary conditions
\begin{eqnarray}\label{eqn:boundconds_lim_tilde}
\displaystyle\lim_{s\rightarrow\pm\infty} \widetilde{u}(s,t) = \widetilde{x}^{\pm}\in\mathcal{P}(\widetilde{H})
\end{eqnarray}
can be seen as connecting orbits between $[\widetilde{x}^+,\widetilde{v}^+]$ and $[\widetilde{x}^-,\widetilde{v}^-]$ in $\widetilde{\mathcal{L}}\widetilde{M}$, where
\begin{eqnarray}\label{eqn:tranfcap}
\widetilde{v}^-\# \widetilde{u}=\widetilde{v}^+.
\end{eqnarray}
The energy of a connecting orbit in this space is given by
\begin{eqnarray}\label{eqn:energy}
E(\widetilde{u})=\displaystyle\int_{-\infty}^{\infty}\displaystyle\int_0^1 ||\partial_s \widetilde{u}||_{\widetilde{M}}^2 dt ds= \mathcal{A}_{\widetilde{H}}([\widetilde{x}^-,\widetilde{v}^-])-\mathcal{A}_{\widetilde{H}}([\widetilde{x}^+,\widetilde{v}^+])
\end{eqnarray}
where $||\partial_s \widetilde{u}||_{\widetilde{M}}$ is $||\partial_s (\pi\circ\widetilde{u})||$ where $||\cdot||$ stands for the norm with respect to $\left<\cdot,\cdot \right>=\omega(\cdot,J\cdot)$. We denote by $\mathcal{M}([\widetilde{x}^-,\widetilde{v}^-],[\widetilde{x}^+,\widetilde{v}^+])$ the space of solutions $\widetilde{u}$ of (\ref{eqn:connectorbit_eq_tilde}-\ref{eqn:boundconds_lim_tilde}) that transform the cappings as in (\ref{eqn:tranfcap}) and that have finite energy $E(\widetilde{u})$.

The Conley-Zehnder index $\MUCZ$ of a non-degenerate periodic solution $[\widetilde{x},\widetilde{v}]\in \widetilde{\mathcal{L}}\widetilde{M}$ satisfies
\begin{eqnarray}\label{eqn:mi_czi}
0\not =|\D_{\widetilde{H}}([\widetilde{x},\widetilde{v}])-\MUCZ([\widetilde{x},\widetilde{v}])|< n.
\end{eqnarray}
and is $\Gamma_1$-invariant (where $\Gamma_1$ is $\pi_1(M) / \ker I_{\theta}$), i.e. 
$$\MUCZ([\widetilde{x},\widetilde{v}])=\MUCZ(a\cdot[\widetilde{x},\widetilde{v}]) \quad\text{for any}\quad a \in \Gamma_1.$$ This index satisfies the following identities:
\begin{itemize}
\item[i)] $\MUCZ([\widetilde{x},\widetilde{v}\#A])=\MUCZ([\widetilde{x},\widetilde{v}])+I_{c_1}(A)$

\item[ii)] dim $\mathcal{M}([\widetilde{x}^-,\widetilde{v}^-],[\widetilde{x}^+,\widetilde{v}^+])= \MUCZ([\widetilde{x}^-,\widetilde{v}^- \# A]) - \MUCZ([\widetilde{x}^+,\widetilde{v}^+ \# A])$
\end{itemize}
for $A\in \pi_2(M)$.

Denote by $\widetilde{\mathcal{P}}_k(\widetilde{H})$ the subset of $\widetilde{\mathcal{P}}(\widetilde{H})$ of periodic solutions with
${\MUCZ([\widetilde{x},\widetilde{v}])=k}$. Consider the chain complex whose $k$-th chain group $C_k(\widetilde{H})$ consists of all formal sums
\begin{eqnarray*}
\sum \xi_{[\widetilde{x},\widetilde{v}]}\cdot [\widetilde{x},\widetilde{v}]
\end{eqnarray*}
with $[\widetilde{x},\widetilde{v}]\in \widetilde{\mathcal{P}}_k(\widetilde{H}),\;\xi_{[\widetilde{x},\widetilde{v}]}\in \mathbb{Z}_2$ and such that, for all $c\in \reals$,
$$
\#\big{\{}[\widetilde{x},\widetilde{v}]\;|\;\xi_{[\widetilde{x},\widetilde{v}]}\not=0,
\;\mathcal{A}_{\widetilde{H}}([\widetilde{x},\widetilde{v}])>c\big{\}}<\infty.
$$
For a generator $[\widetilde{x},\widetilde{v}]$ in $C_k(\widetilde{H})$, the boundary operator $\partial_k$ is defined as follows
$$
\partial_k([\widetilde{x},\widetilde{v}])=\displaystyle\sum_{\MUCZ([\widetilde{y},\widetilde{w}])=k-1} n_2([\widetilde{x},\widetilde{v}],[\widetilde{y},\widetilde{w}])[\widetilde{y},\widetilde{w}]
$$
where $n_2([\widetilde{x},\widetilde{v}],[\widetilde{y},\widetilde{w}])\in \mathbb{Z}_2$ is the modulo-$2$ reduction of the number of elements in the quotient space $\mathcal{M}([\widetilde{x},\widetilde{v}],[\widetilde{y},\widetilde{w}])/\reals$. The boundary operator $\partial$ satisfies $\partial^2=0$ and we have the Floer-Novikov homology groups
\begin{eqnarray}
HFN_k(\theta_t)=\frac{\ker\partial_k}{\im \partial_{k+1}}.
\end{eqnarray}
Moreover, this homology is invariant under exact deformations of the closed form $\theta_t$ (see \cite[Theorem~4.3]{LO_fixedpts95}) and hence two paths with the same flux have isomorphic associated Floer-Novikov homology groups.

\subsubsection{Filtered Floer-Novikov Homology}
The (total) chain Floer complex $C_*(\widetilde{H})=\colon$ $C_*^{(-\infty, \infty)}(\widetilde{H})$ admits a filtration by $\reals$. Define $\CS(\widetilde{H})$ the set of critical values of the functional $\mathcal{A}_{\widetilde{H}}$ (defined in (\ref{eqn:action_flnov})) which is called the \emph{action spectrum} of $\{\phi_t\}$. For each
$b\in(-\infty,\infty]$ outside $\CS(\widetilde{H})$, the chain complex $C_*^{(-\infty, b)}(\widetilde{H})$ is generated by equivalence classes of capped loops $[(\widetilde{x}, \widetilde{v})]$ with action $\A_{\widetilde{H}}$ less than~$b.$ For $-\infty\leq a < b \leq \infty$ outside $\CS(\widetilde{H}),$ set
$$
C_*^{(a,b)}(\widetilde{H}) :=C_*^{(-\infty,b)}(\widetilde{H}) / C_*^{(-\infty,a)}(\widetilde{H}).
$$
The boundary operator $\p\colon C_*(\widetilde{H})\rightarrow C_{*-1}(\widetilde{H})$ descends to $C_*^{(a,b)}(\widetilde{H})$ and hence the \emph{filtered Floer-Novikov homology} $HFN_*^{(a,b)}(\theta_t)$ is well defined.

This construction also extends by continuity to all symplectomorphisms in $\mbox{Symp}_0(M,\omega).$ For a path $\{\phi_t\}$ connecting an arbitrary $\phi\in\mbox{Symp}_0(M,\omega)$ to the identity, set 
\begin{eqnarray}\label{eqn:homol_deg}
HFN_*^{(a,b)}(\theta_t)\colon=
HFN_*^{(a,b)}(\theta_t^{\prime})
\end{eqnarray}
where $\theta_t$ is the one-form (defined as in Lemma~\ref{lemma:deformation}) associated with the path $\{\phi_t\}$ and $\theta_t^{\prime}$ is the one-form associated with a non-degenerate perturbation $\{\phi^{\prime}_t\}$ of $\{\phi_t\}$. Here $-\infty\leq a<b\leq \infty$ are outside the action spectrum of $\{\phi_t\}$. 
Observe that since the symplectic manifold $(M,\omega)$ and the flux are rational (in the sense of Sections~\ref{section:sympmnfd}~and~\ref{section:sympmrpfisms}), the action spectrum is nowhere dense (see e.g. \cite{HZ11, Schw00}).

\subsection{Quantum homology} The quantum homology of $M$, $HQ_*(M)$, is an algebra over the Novikov ring, $\Lambda$. In this section we recall their definitions; see  \cite[Chapter $11$]{MS12} for more details. Here we follow \cite[Section $2.2$]{GG:hyp12}.

\subsubsection{Novikov ring and Quantum homology}\label{subsubsection:novikov}In the case when $M$ is strictly monotone, the \emph{Novikov ring}, $\Lambda$, is the group algebra of a group $\Gamma$ over $\Z_2$, $\Z_2[\Gamma]$. The group $\Gamma$ is the quotient of $\pi_2(M)$ by the equivalence relation $\sim$ where $A\sim B$ if $I_{c_1}(A)=I_{c_1}(B)$, or equivalently, if  $I_{\omega}(A)=I_{\omega}(B)$, i.e.
$$
\Gamma=\frac{\pi_2(M)}{\ker I_{\omega}}=\frac{\pi_2(M)}{\ker I_{c_1}}.
$$
An element in $\Lambda$ is a formal finite linear combination,
$$
\sum \alpha_{\small A} e^A,
$$
where $\alpha_{\small A} \in \Z_2$. The degree of $e^A$, for $A\in \Gamma$, is $I_{c_1}(A)$ which grades the ring $\Lambda$. We have $\Gamma \simeq \Z$ and denote by $A_0$ the generator of $\Gamma$ with $I_{c_1}(A_0)=-2c_1^{\min}$. Then $\mathfrak{q}:=e^{A_0}\in \Lambda$ has degree $-2c_1^{\min}$ and the Novikov ring is the ring of Laurent polynomials $\Z_2[\mathfrak{q}^{-1},\mathfrak{q}]$.\\

The \emph{quantum homology} of $M$ is defined by
$$
HQ_*(M)= H_*(M)\otimes\Lambda
$$
(where $\Lambda$ is the Novikov ring) where the degree of an element ${\alpha\otimes e^A}$, with $\alpha\in H_*(M)$ and $A\in\Gamma$, is ${\deg(\alpha)+I_{c_1}(A)}$. The product structure is given by the \emph{quantum product}. When $\alpha$ and $\beta$ are ordinary homology classes the quantum product is 
\begin{eqnarray}\label{eqn:quantum_product}
\alpha*\beta=\displaystyle\sum_{A\in\Gamma}(\alpha*\beta)_{\small A} e^A,
\end{eqnarray}
where $(\alpha*\beta)_{\small A}\in H_*(M)$ is defined using Gromov-Witten invariants of $M$ and its degree is ${\deg(\alpha)+\deg(\beta)-2n-I_{c_1}(A)}$. Then $${\deg(\alpha*\beta)=\deg(\alpha)+\deg(\beta)-2n}.$$ When $A=0$, $(\alpha*\beta)_0=\alpha\cap\beta$, where $\cap$ is the intersection product of ordinary homology classes.

Recall that the summation in (\ref{eqn:quantum_product}) can be restricted to the cone $I_{\omega}(A)\leq 0$ and, under the hypotheses on $M$, we can rewrite (\ref{eqn:quantum_product}) as 
\begin{eqnarray}\label{def:quantumprod}
\alpha*\beta=\alpha\cap\beta+\displaystyle\sum_{k>0}(\alpha*\beta)_k \mathfrak{q}^k,
\end{eqnarray}
where the sum is finite and the degree of $(\alpha*\beta)_k$ is $\deg((\alpha*\beta)_k)=\deg(\alpha)+\deg(\beta)-2n+2c_1^{\min}k$.

The unit in the algebra $HQ_*(M)$ is the fundamental class $[M]$ and, for $a\in \Lambda$ and $\alpha\in H_*(M)$, 
$$a\alpha=(a[M])*\alpha$$ where degree of $a\alpha$ is $\deg(a\alpha)=\deg(a)+\deg(\alpha)$. Then the ordinary homology $H_*(M)$ is canonically embedded in $HQ_*(M)$.

The map $I_{\omega}$ (initially defined on $\pi_2(M)$ in (\ref{eqn:I_omega})) can be defined on $HQ_*(M)$ by
$$
I_{\omega}(\alpha)=\max\{I_{\omega}(A)\;|\; \alpha_{\small A} \not= 0 \}=\max\{-h_0k\;|\;\alpha_k\not=0\}
$$
where $\alpha=\sum\alpha_{\small A} e^A=\sum\alpha_k \mathfrak{q}^k$ and it satisfies
\begin{eqnarray*}
I_{\omega}(\alpha+\beta)\leq\max\{I_{\omega}(\alpha),I_{\omega}(\beta)\}
\end{eqnarray*}
and
\begin{eqnarray}\label{eqn:I_omega_inequality_2}
I_{\omega}(\beta*\alpha)\leq I_{\omega}(\alpha)+I_{\omega}(\beta).
\end{eqnarray}

\begin{Example}\label{example:quantumprod}
Consider $M=\CP^n,\;\alpha=[\CP^{n-1}]\in H_{2n-2}(M)$ and $\beta=[\text{pt}]\in H_{0}(M)$ (where $[\text{pt}]$ is the class of a point in $\CP^n$). The fact that there is a unique line through any two points is reflected in the identity $(\beta*\alpha)_A=[\CP^n]$ where $c_1(A)=c_1^{\min}$. Hence, in $HQ_*(\CP^n)$, $\beta*\alpha=\mathfrak{q}[\CP^n].$ 

Similarlly, for    
$M= \CP^n\times P^{2m}$, where $P^{2m}$ is symplectically aspherical,
$$
\alpha=[\CP^{n-1}\times P^{2m}]\in H_{2n+2m-2}(M)
$$ 
and 
$$
\beta=[\text{pt} \times P^{2m}]\in H_{2m}(M),
$$ 
the quantum product $\beta*\alpha$ satisfies the homological condition (\ref{l:1}). (See \cite{MS12} for these and further computations.)
%Then the only term which appears in (\ref{def:quantumprod}) corresponds to $k=1$ and $\deg((\beta*\alpha)_1)=(2n+2m-2)+2m-2(n+m)+2(n+1)=2(n+m)=\dim M$.
%Hence, $\beta*\alpha$ is a multiple of $[M]$ where the multiple is given by the invariant $\Phi_{A_0} (\beta, \alpha, [\text{pt}^{'}]\otimes[\text{pt}^{'}])$ defined in \cite{MS12}, where $[\text{pt}]$ is the class of a point in $P^{2m}$. The fact that there is a unique line through any two points is reflected in the identity
%$$
%\Phi_{A_0} ([\CP^{n-1}]\otimes [P^{2m}], [\text{pt}]\otimes\otimes [P^{2m}], [\text{pt}]\otimes [\text{pt}^{'}])=1
%%$$
%(see \cite{MS12} for these computations) and then the homological assumption in Theorem~\ref{maintheo}
%$$
%\beta*\alpha=\mathfrak{q}[M]
%$$
%holds.
\qed
\end{Example} 
\subsubsection{Quantum Product Action}\label{section:capproduct} We describe an action of the quantum homology on the filtered Floer-Novikov homology. We follow  \cite[Section~2.3]{GG:hyp12} for the Floer-Novikov setting; see \cite[Section~3]{LO_cup96} for more details. 

Let $\phi$ be a non-degenerate symplectomorphism, $J$ be a generic almost complex structure and $[\sigma]$ be a class in $H_*(M)$. Denote by $\mathcal{M}([\widetilde{x},\widetilde{v}],[\widetilde{y},\widetilde{w}];\sigma)$ the moduli space of solutions 
$\widetilde{u}$ of (\ref{eqn:connectorbit_eq_tilde}-\ref{eqn:boundconds_lim_tilde}) that transform the cappings as in (\ref{eqn:tranfcap}) with ${[\widetilde{x},\widetilde{v}], [\widetilde{y},\widetilde{w}] \in \widetilde{\mathcal{P}}(\widetilde{H})}$ and such that $u(0,0)\in \sigma$ where $\sigma$ is a generic cycle representing $[\sigma]$ and $\pi \circ\widetilde{u}=u$.

Then the dimension of this moduli space is
$$
\mbox{dim} \;\mathcal{M} ([\widetilde{x},\widetilde{v}],[\widetilde{y},\widetilde{w}];\sigma)=\MUCZ([\widetilde{x},\widetilde{v}])-
        \MUCZ([\widetilde{y},\widetilde{w}])
-\mbox{codim}(\sigma).
$$
Let $m([\widetilde{x},\widetilde{v}],[\widetilde{y},\widetilde{w}];\sigma)\in\Z_2$ be $\# \mathcal{M}([\widetilde{x},\widetilde{v}],[\widetilde{y},\widetilde{w}];\sigma)$ mudulo 2 when this moduli space is zero-dimensional and zero otherwise. For any $c,\;c'\not\in S(\widetilde{H})$, there is a map 
$$
\Phi_{\sigma}:C_*^{(c,c')}(\widetilde{H})\rightarrow C_{*-\mbox{codim}(\sigma)}^{(c,c')}(\widetilde{H})
$$
induced by
$$
\Phi_{\sigma}([\widetilde{x},\widetilde{v}])=\displaystyle\sum_{[\widetilde{y},\widetilde{w}]} m([\widetilde{x},\widetilde{v}], [\widetilde{y},\widetilde{w}]; \sigma) [\widetilde{y},\widetilde{w}].
$$
This map commutes with the Floer-Novikov differential $\p$ and descends to a map
$$
\Phi_{[\sigma]}:HFN_*^{(c,c')}(\theta_t)\rightarrow HFN_{*-\mbox{codim}(\sigma)}^{(c,c')}(\theta_t)
$$
(this is independent of the cycle which represents $[\sigma]$).
The map
\begin{eqnarray}\label{eqn:action_induced}
\Phi_{\mathfrak{q}^{\nu}\sigma}([\widetilde{x},\widetilde{v}]):=\displaystyle\sum_{[\widetilde{y},\widetilde{w}]}m(\mathfrak{q}^{\nu}[\widetilde{x},\widetilde{v}], [\widetilde{y},\widetilde{w}]; \sigma)[\widetilde{y},\widetilde{w}]
\end{eqnarray}
induces an action of $\alpha=\mathfrak{q}^{\nu}[\sigma]\in HQ_*(M)$ 
\begin{eqnarray}\label{eqn:action_quantum}
\Phi_{\alpha}: HFN_*^{(c,c')}(\theta_t)\rightarrow HFN_{*-2n+\deg(\alpha)}^{(c,c')+I_{\omega}(\alpha)}(\theta_t).
\end{eqnarray}

In (\ref{eqn:action_induced}), $\mathfrak{q}$ is as in Section~\ref{subsubsection:novikov} and $\mathfrak{q}^{\nu}[\widetilde{x},\widetilde{v}]$ is the element $[\widetilde{x},\widetilde{u}]$ in $ \widetilde{\mathcal{P}}(\widetilde{H})$ where $\widetilde{u}$ is 
obtained by ``recapping $\widetilde{v}$" in the following way: $u=v\# (\nu A_0)$. Here $A_0$ is the generator of the group $\Gamma$ defined in Section~\ref{subsubsection:novikov}, $u=\pi\circ\widetilde{u}$ and $v=\pi\circ\widetilde{v}$.

The map $\Phi_{\alpha}$ can be extended to all $\alpha\in HQ_*(M)$, by linearity over $\Lambda$, so that (\ref{eqn:action_quantum}) holds. 

The maps $\Phi_{\alpha}$, for all $\alpha\in HQ_*(M)$, also give an action of the quantum homology on the filtered Floer-Novikov homology. 

The action has the following properties:
$$
\Phi_{[M]}=id
$$
and 
\begin{eqnarray}\label{eqn:multip_quantum}
\Phi_{\beta}\Phi_{\alpha}=\Phi_{\beta*\alpha}.
\end{eqnarray}

\begin{Remark}
Observe that in the multiplicativity property (\ref{eqn:multip_quantum}), in general, the maps on the two sides of the equality have different target spaces. We should understand the identity in (\ref{eqn:multip_quantum}) as that the following diagram commutes for any interval $(d,d')$
\begin{equation}\label{diagram:cap}
\xymatrix{
{HFN_*^{(c,c^{\prime})}(\theta_t)} \ar[r]^{\Phi_{\alpha}}\ar[dr]_{\Phi_{\beta*\alpha}}
            & {HFN_{*-2n+\deg(\alpha)}^{(c,c^{\prime})+I_{\omega}(\alpha)}(\theta_t)} \ar[r]^{\Phi_{\beta}} & {HFN_{*-4n+\deg(\alpha)+\deg(\beta)}^{(c,c^{\prime})+I_{\omega}(\alpha)+I_{\omega}(\beta)}}(\theta_t) \ar[d] \\
{}                     & {HFN_{*-2n+\deg(\beta*\alpha)}^{(c,c^{\prime})+I_{\omega}(\beta*\alpha)}(\theta_t)}             \ar[r]   &  {HFN_{*-2n+\deg(\beta*\alpha)}^{(d,d')}(\theta_t)}
}
\end{equation}
where $d\geq c+I_{\omega}(\alpha)+I_{\omega}(\beta)$ and $d'\geq c'+I_{\omega}(\alpha)+I_{\omega}(\beta)$.  The vertical arrow and the bottom horizontal arrow are the natural quotient-inclusion maps whose existence is guaranteed by the choice of $d$ and $d'$ and (\ref{eqn:I_omega_inequality_2}).\qed
\end{Remark}

\section{Proof of the Main Result}\label{section:proof}
\subsection{Ball-Crossing Energy Theorem}\label{section:BCEthm} Here, we describe the key property of hyperbolic periodic orbits which supports the proof of the main theorem (see \cite[Section~3]{GG:hyp12} for more details including the proof of the Ball-crossing energy Theorem).

Let $\phi$ be a symplectomorphism isotopic to the identity on a symplectic manifold $(M,\omega)$ and fix a one-periodic in time almost complex structure $J$ compatible with $\omega$. For a closed domain $\Sigma \subset \reals \times S^1_k$ (i.e. a closed subset with non-empty interior), where $S^1_k:=\reals / k\Z$, the energy of a solution ${\widetilde{u}:\Sigma \rightarrow \widetilde{M}}$ of the equation (\ref{eqn:connectorbit_eq_tilde}) is, by definition, 
$$
E(\widetilde{u}):= \int_{\Sigma} ||\partial_s \widetilde{u}||^2_{\widetilde{M}}\; dt ds
$$
where $||\partial_s \widetilde{u}||_{\widetilde{M}}$ is defined as in (\ref{eqn:energy}).

Let $\gamma$ be a hyperbolic one-periodic solution of (\ref{eqn:de_theta}) in $M$ and $\widetilde{\gamma}$ a lift of $\gamma$ to $\widetilde{M}$, then $\widetilde{\gamma}\in \mathcal{P}(\widetilde{H})$ hyperbolic. Recall the definition of the covering space $\widetilde{M}$ in Section~\ref{section:tailedcapped} and of the Hamiltonian $\widetilde{H}$ associated with $\phi$ in~(\ref{eqn:H_tilde}).

A solution ${\widetilde{u}:\Sigma \rightarrow \widetilde{M}}$ of the equation (\ref{eqn:connectorbit_eq_tilde}) is said to be asymptotic to $\widetilde{\gamma}^k$ as $s\rightarrow \infty$ if $\Sigma$ contains a cylinder $[s_0,\infty)\times S^1_k$ and $\widetilde{u}(s,t)\rightarrow \widetilde{\gamma}^k(t)$ $C^{\infty}$-uniformly in $t$ as $s\rightarrow \infty$.

Consider a \emph{sufficiently small} closed neighborhood U of $\gamma$ with smooth boundary and define $\widetilde{U}:=\pi^{-1}(U)$.

\begin{Theorem}\cite[Ball-Crossing Energy Theorem]{GG:hyp12}
\label{thm:ballcrossing}
\

There exists a constant $c_{\infty}>0$ (independent of $k$ and $\Sigma$) such that for any solution $\widetilde{u}$ of the equation (\ref{eqn:connectorbit_eq_tilde}), with $\widetilde{u}(\partial\Sigma)\subset\partial \widetilde{U}$ and $\partial\Sigma\not=\emptyset$, which is asymptotic to $\widetilde{\gamma}^k$ as $s\rightarrow\infty$, we have
\begin{eqnarray}\label{eqn:energy_bound}
E(\widetilde{u})>c_{\infty}.
\end{eqnarray}
Moreover, the constant $c_{\infty}$ can be chosen to make (\ref{eqn:energy_bound}) hold for all $k$-periodic almost complex structures (with varying $k$) $C^{\infty}$-close to $\widetilde{J}$ uniformly on $\reals\times \widetilde{U}$.
\end{Theorem}

\subsection{Proof of Theorem~\ref{maintheo}}\label{proof:maintheo}
As mentioned in the introduction, to simplify the setting, we may assume that $\gamma$ is a constant orbit.  This is due to the fact that there exists a one periodic loop of Hamiltonian diffeomorphisms $\psi^t$ defined on a neighborhood of $\gamma$ such that $\psi^t(\gamma(0))=\gamma(t)$. We may think of $\gamma(t)\equiv \gamma(0)$ as a fixed point of $(\psi^t)^{-1}\circ\phi_t$ (see \cite[Section~5.1]{Gi:conley} for more details).

Denote by $\widetilde{H}$ the Hamiltonian associated with the symplectic isotopy $\phi_t$ in the sense of (\ref{eqn:H_tilde}). Consider, by (\ref{eqn:mean_homo}), an iteration $l$ of $\widetilde{H}$ so that the mean index of the $l$-th iteration of $\gamma$ is divisible by $2c_1^{\min}$ (independently of the attached capping). Then, by (\ref{eqn:mean_recap}), the mean index of this iteration of $\gamma$ is zero for some capping. Add a constant to $\widetilde{H}^{\natural l}$ so that the action of a lift $\widehat{\gamma}$ of the $l$-th iteration of $\gamma$ with the above capping is zero. We write $\widetilde{H}$ and $\gamma$ for the obtained Hamiltonian and the $l$-th iteration of $\gamma$, respectively, and we have: 
\begin{eqnarray}\label{eqn:unif_gamma}
\D_{\widetilde{H}}({\widehat{\gamma}})=0=\mathcal{A}_{\widetilde{H}}({\widehat{\gamma}}).
\end{eqnarray}

\begin{Remark}
%\
%\begin{enumerate}
%\item To ensure condition (\ref{eqn:unif_gamma}), we may have to consider an iteration of $\phi$ which we continue denoting by $\phi$: in fact, we pass, if necessary, to an iteration so that the mean index of $\gamma$ is divisable by $2c_1^{\min}$ (independently of the attached capping). Then the mean index is zero for some capping and finally by adding a constant to the obtained associated Hamiltonian we can assume the action is also zero. 
%\item 
Observe that since $\gamma$ is hyperbolic, the mean index $\D_{\widetilde{H}}({\widehat{\gamma}})$ is equal to the Conley-Zehnder index $\MUCZ({\widehat{\gamma}})$ and hence $\MUCZ({\widehat{\gamma}})=0$.\qed
%\end{enumerate}
\end{Remark}

We reason by contradiction and suppose $\phi$ has finitely many periodic points. By our assumptions we have 
$$
\left<\widetilde{\text{Flux}}([\phi_t]),\pi_1(M)\right>=\frac{p}{q}h_0\Z,
$$
where $p$ and $q$ are relatively prime natural numbers (and $q$ is non-zero). Notice that, for all $k\in\N$, we have 
$$
\left<\widetilde{\text{Flux}}([\phi_t^{k}]),\pi_1(M)\right>=\frac{kp}{q}h_0\Z.
$$
\begin{Remark}\label{rmk:flux_it}
Observe that $\widetilde{\mbox{Flux}}([\phi_t^k])=k\widetilde{\mbox{Flux}}([\phi_t])$.\qed
\end{Remark}
We fix a prime number $k$ large enough so that 
\begin{eqnarray}\label{eqn:p}
(2kp-3)h_0 - \lambda(n+1)>0
\end{eqnarray}
and from now on we simply write $\phi_t$ and $p$ for $\phi_t^k$ and $kp$, respectively. 

The $q$-th iteration $\phi^q$ also has finitely many periodic points and we denote them by $x_1,\ldots,x_m$. We write $\phi_t^q$ for $\phi_t$ and (up to further replacing of $\phi^k_t$ by $\phi_t$) the curves $x_1,\ldots,x_m$ can be assumed to be one-periodic.

%Consider an iteration of $\phi$, still denoted by $\phi$, so that $p$ is sufficiently large, namely,
%\begin{eqnarray}\label{eqn:p}
%(2p-3)h_0 - \lambda(n+1)>0
%\end{eqnarray}
%where $\lambda$ is the monotonicity constant of $M$. The $r$-th iteration $\phi^r$ (where $r$ is defined in Theorem~\ref{maintheo}) has finitely many periodic points and we denote them by $x_1,\ldots,x_m$.

%\begin{Remark}\label{rmk:flux_it}
%\
%\begin{enumerate}
%\item Observe that $\widetilde{\mbox{Flux}}(\{\phi_t^r\})=r\widetilde{\mbox{Flux}}(\{\phi_t\})$.
%\item The periodic orbit $\gamma^r$ of $\phi^r$ is hyperbolic and we keep the notation $\gamma$ for this orbit and $\phi$ for the iteration $\phi^r$.
%\end{enumerate}
%\end{Remark}

Consider an almost complex structure $J^{'}$ and denote by $\widetilde{J^{'}}$ the corresponding almost complex structure on $\widetilde{M}$.  Let $U$ be a closed neighborhood of $\gamma$ such that $\phi$ has no periodic orbit intersecting $U$ except $\gamma$. By Theorem~\ref{thm:ballcrossing}, there exists a constant $c_{\infty}>0$ such that, for all $k$, the energy of any solution of (\ref{eqn:connectorbit_eq_tilde}) of period $k$ asymptotic to $\widetilde{\gamma}^k$ as $s\rightarrow\infty$ is greater than $c_{\infty}$.

For each $i=1,\ldots,m$, attach a capping $v_i$ to the loop $x_i$, fix a lift $[(\widetilde{x}_i,\widetilde{v}_i)]=:\widehat{x}_i\in \LL\MM$ and define
\begin{itemize}
\item[] $a_i:= \mathcal{A}_{\widetilde{H}}(\widehat{x}_i)$ (mod $h_0$) in $ S^1_{h_0}$ 
\item[] $\underline{a_i} :=\underline{\mathcal{A}}_{\widetilde{H}}(\widehat{x}_i)$ (mod $p h_0$) in $ S^1_{p h_0}$.
\end{itemize}

\begin{Remark}\label{rmk:ind_action}
Observe that $a_i$ and $\underline{a_i}$ are independent of the initially attached capping and fixed lift. (The second follows from Remark~\ref{rmk:flux_it}.)\qed
\end{Remark}

Take $\epsilon,\;\delta>0$ small, namely,
\begin{eqnarray}\label{epsdelt}
2(\epsilon+\delta)<\lambda \mbox{ and } \epsilon<c_{\infty}.
\end{eqnarray}
Then, by Kronecker's Theorem, there exists $k$ (arbitrarily large) such that for all $i=1,\ldots,m$
\begin{eqnarray}\label{eqn:kronecker}
||k a_i||_{h_0}<\epsilon\quad \mbox{and} \quad ||k\underline{a_i}||_{p h_0}<\delta.
\end{eqnarray}
We denote by $||a||_{h}\in [0,h / 2]$ the distance from $a\in S^1_{h}=\reals/{h\Z}$ to $0$. Observe that $k$ depends on $\epsilon$ (and $\delta$), hence it depends on $c_{\infty}$ and also on the neighborhood $U$.

%Consider a non-degenerate perturbation $\phi^{'}$ of $\phi^k$ such that (\ref{eqn:homol_deg}) holds and the Hamiltonian $\widetilde{K}$ associated with $\phi^{'}$ (in the sense of (\ref{eqn:H_tilde})) satisfies the following properties:
Consider a non-degenerate perturbation $\phi_t^{'}$ of $\phi_t^{k}$ such that (\ref{eqn:homol_deg}) holds and such that the Hamiltonian $\widetilde{K}$ associated with $\phi_t^{'}$ as in (\ref{eqn:H_tilde}) satisfies the following properties:

\begin{enumerate}\label{enumerate:K}
\item \label{enumerate:K2}$\widetilde{K}$ coincides with $\widetilde{H}^{\natural k}$ on the neighborhood $U$,
\item \label{enumerate:K3}$\widetilde{K}$ is $k$-periodic and non-degenerate and
\item \label{enumerate:K1}$\widetilde{K}$ is sufficiently $C^2$-close to $\widetilde{H}^{\natural k}$.
\end{enumerate}
Here, we consider $\widetilde{K}$ sufficiently $C^2$-close to $\widetilde{H}^{\natural k}$ in order to have the existence of $k$ (large) such that for all $\widehat{x}$ $k$-periodic solution of $\widetilde{K}$ 
\begin{eqnarray}\label{eqn:action_epsbounds}
\big\|\mathcal{A}^{h_0}_{\widetilde{K}}(\widehat{x})\big\|_{h_0}<\epsilon
\end{eqnarray}
and
\begin{eqnarray}\label{eqn:augmented_bounds}
\mbox{either } \big|\underline{\mathcal{A}}_{\widetilde{K}}(\widehat{x})\big|<\delta
\mbox{ or } \big|\underline{\mathcal{A}}_{\widetilde{K}}(\widehat{x})\big|>(p-1)h_0
\end{eqnarray}
where $\mathcal{A}^{h_0}_{\widetilde{K}}(\widehat{x})$ stands for $\mathcal{A}_{\widetilde{K}}(\widehat{x})$ mod $h_0$. Note that, as long as $\delta<h_0$, conditions (\ref{eqn:augmented_bounds}) and (\ref{eqn:action_epsbounds}) follow from (\ref{eqn:kronecker}). Observe that if $\phi^k$ is non-degenerate, we can take $\phi^{'}=\phi^k$. 

For any $k$-periodic almost complex structure $\widetilde{J}$ sufficiently close to (the $k$-periodic extension of) $\widetilde{J^{'}}$, all non-trivial solutions $\widetilde{u}\colon \reals \times S^1_k \rightarrow \widetilde{M}$ of the equation (\ref{eqn:connectorbit_eq_tilde}) for the pair ($\phi^{'}$, $\widetilde{J}$) asymptotic to $\widetilde{\gamma}^k$ as $s\rightarrow\infty$ have energy greater than $c_{\infty}$.

%all non-trivial $k$-periodic solutions of the equation (\ref{eqn:connectorbit_eq_tilde}) for the pair ($\phi^{'}$, $\widetilde{J}$) asymptotic to $\widetilde{\gamma}^k$ as $s\rightarrow\infty$ have energy greater than $c_{\infty}$.

\begin{Lemma}\emph{{\cite[Lemma~4.1]{GG:hyp12}}}\label{lemma:tau}
Let $\tau := (p-1)h_0 - \frac{\lambda}{2}(n+1)$. The orbit ${\widehat{\gamma}}^k$ is not connected by a solution of (\ref{eqn:connectorbit_eq_tilde}) to any $\widehat{x}\in \widetilde{\mathcal{P}}(\widetilde{K})$  with Conley-Zehnder index $\pm 1$ and action in $(-\tau,\tau)$.

In particular, ${\widehat{\gamma}}^k$ is closed in $C^{(-\tau,\tau)}_*(\widetilde{K})$ and $0\not=[{\widehat{\gamma}}^k]\in HFN^{(-\tau,\tau)}_*(\theta^{\prime}_t)$. Moreover, ${\widehat{\gamma}}^k$ must enter every cycle representing its homology class $[{\widehat{\gamma}}^k]$ in $HFN^{(-\tau,\tau)}_*(\theta^{\prime}_t)$.
\end{Lemma}
\proof
Assume the orbit ${\widehat{\gamma}}^k$ is connected, by a solution $\widetilde{u}$ of (\ref{eqn:connectorbit_eq_tilde}), to some $\widehat{x}\in \widetilde{P}(\widetilde{K})$ with index $\MUCZ(\widehat{x})=\pm 1$ with action in $(-\tau,\tau)$.

Consider the first case in (\ref{eqn:augmented_bounds}), i.e. $\big|\underline{\mathcal{A}}_{\widetilde{K}}(\widehat{x})\big|<\delta$: since
\begin{itemize}
\item[i)] $\big\|\mathcal{A}^{h_0}_{\widetilde{K}}(\widehat{x})\big\|_{h_0}<\epsilon$ (by (\ref{eqn:action_epsbounds})),
\item[ii)]$E(\widetilde{u})>c_{\infty}>\epsilon$ (by Theorem~\ref{thm:ballcrossing} and (\ref{epsdelt})) and
\item[iii)]$\mathcal{A}_{\widetilde{K}}({\widehat{\gamma}}^k)=0$ (by (\ref{eqn:unif_gamma})),
\end{itemize}
we have
$$
\big|\mathcal{A}_{\widetilde{K}}(\widehat{x})\big|>h_0 - \epsilon.
$$
Then, by the definition of augmented action (\ref{eqn:augmentedaction}) and since
\begin{itemize}
\item[i)] $\big|\underline{\mathcal{A}}_{\widetilde{K}}(\widehat{x})\big|<\delta$ and
\item[ii)] $2(\epsilon +\delta)<\lambda$ (by \ref{epsdelt}),
\end{itemize}
we have
$$
\big|\D_{\widetilde{K}}(\widehat{x})\big|> \frac{2}{\lambda}(h_0-\epsilon-\delta)=2c_1^{\min}-\frac{2(\epsilon+\delta)}{\lambda} >2c_1^{\min}-1.
$$
Thus, by (\ref{eqn:mi_czi}),
$$
\big|\MUCZ(\widehat{x})\big|>2c_1^{\min}-1-n\geq n+2-1-n=1
$$
where the second inequality follows from the requirement that $c_1^{\min}\geq n/2 +1$. We obtained a contradiction since $\MUCZ(\widehat{x})=\pm 1$.

Consider now the second case in (\ref{eqn:augmented_bounds}), i.e. $\big|\underline{\mathcal{A}}_{\widetilde{K}}(\widehat{x})\big|>(p-1)h_0$: by the definition of augmented action (\ref{eqn:augmentedaction}), we obtain
$$
\big|\mathcal{A}_{\widetilde{K}}(\widehat{x})\big|> (p-1)h_0 - \frac{\lambda}{2} \big|\D_{\widetilde{K}}(\widehat{x})\big| >
(p-1)h_0-\frac{\lambda}{2}(n+1)=:\tau
$$
where the second inequality follows from the fact that $\big|\D_{\widetilde{K}}(\widehat{x})\big|<n+1$ (which holds since $\MUCZ(\widehat{x})=\pm 1$ and by (\ref{eqn:mi_czi})). Hence the action of $\widehat{x}$ is outside the interval $(-\tau,\tau)$ and we obtained a contradiction. \qed\\

The previous lemma also holds for $\mathfrak{q}{\widehat{\gamma}}^k$ with the shifted range of actions $(-\tau,\tau)-h_0.$
For an interval $(a,b)$ containing the interval $[-h_0,0]$ and contained in the intersection of the action intervals $(-\tau,\tau)$ and $(-\tau,\tau)-h_0$,
Lemma~\ref{lemma:tau} holds for both \emph{tailed}-capped orbits ${\widehat{\gamma}}^k$ and $\mathfrak{q}{\widehat{\gamma}}^k$ and the interval $(a,b).$
\begin{Remark}
Observe that the existence of such an interval $(a,b)$ is guaranteed by $-\tau<\tau-h_0$ and our initial assumption on $p$, namely, $(2p-3)h_0 - \lambda(n+1)>0$ (\ref{eqn:p}).\qed
\end{Remark}
For the sake of completeness, we state this result in the following lemma.
\begin{Lemma}\label{lemma:qgamma}
The orbits ${\widehat{\gamma}}^k$ and $\mathfrak{q}{\widehat{\gamma}}^k$ are not connected by a solution of (\ref{eqn:connectorbit_eq_tilde}) to any $\widehat{x}\in \widetilde{\mathcal{P}}(\widetilde{K})$  with Conley-Zehnder index $\pm 1$ and action in an interval $(a,b)$ such that
\begin{eqnarray*}
[-h_0,0]\subset(a,b)\subset (-\tau,\tau)\cap (-\tau-h_0,\tau-h_0).
\end{eqnarray*}

In particular, ${\widehat{\gamma}}^k$ and $\mathfrak{q}{\widehat{\gamma}}^k$ are closed in $C^{(a,b)}_*(\widetilde{K})$ and $[{\widehat{\gamma}}^k]\not=0\not=[\mathfrak{q}{\widehat{\gamma}}^k]\in HFN^{(a,b)}_*(\theta^{\prime}_t)$. Moreover, the orbits ${\widehat{\gamma}}^k$ and $\mathfrak{q}{\widehat{\gamma}}^k$ must enter every cycle representing their homology classes, respectively $[{\widehat{\gamma}}^k]$ and $\mathfrak{q}[{\widehat{\gamma}}^k]$, in $HFN^{(a,b)}_*(\theta^{\prime}_t)$.
\end{Lemma}

Recall that, by (\ref{eqn:action_epsbounds}), all periodic orbits of $\phi^{'}$ have action values in the $\epsilon$-neighborhood of $h_0\Z$. With the following lemma we obtain a contradiction and the main result follows.

\begin{Lemma}\emph{{\cite[Lemma~4.2]{GG:hyp12}}}
The symplectomorphism $\phi^{'}$ has a periodic orbit with action outside the $\epsilon$-neighborhood of $h_0\Z$.
\end{Lemma}

\begin{figure}[htb!]
  \centering
  \def\svgwidth{300pt}
  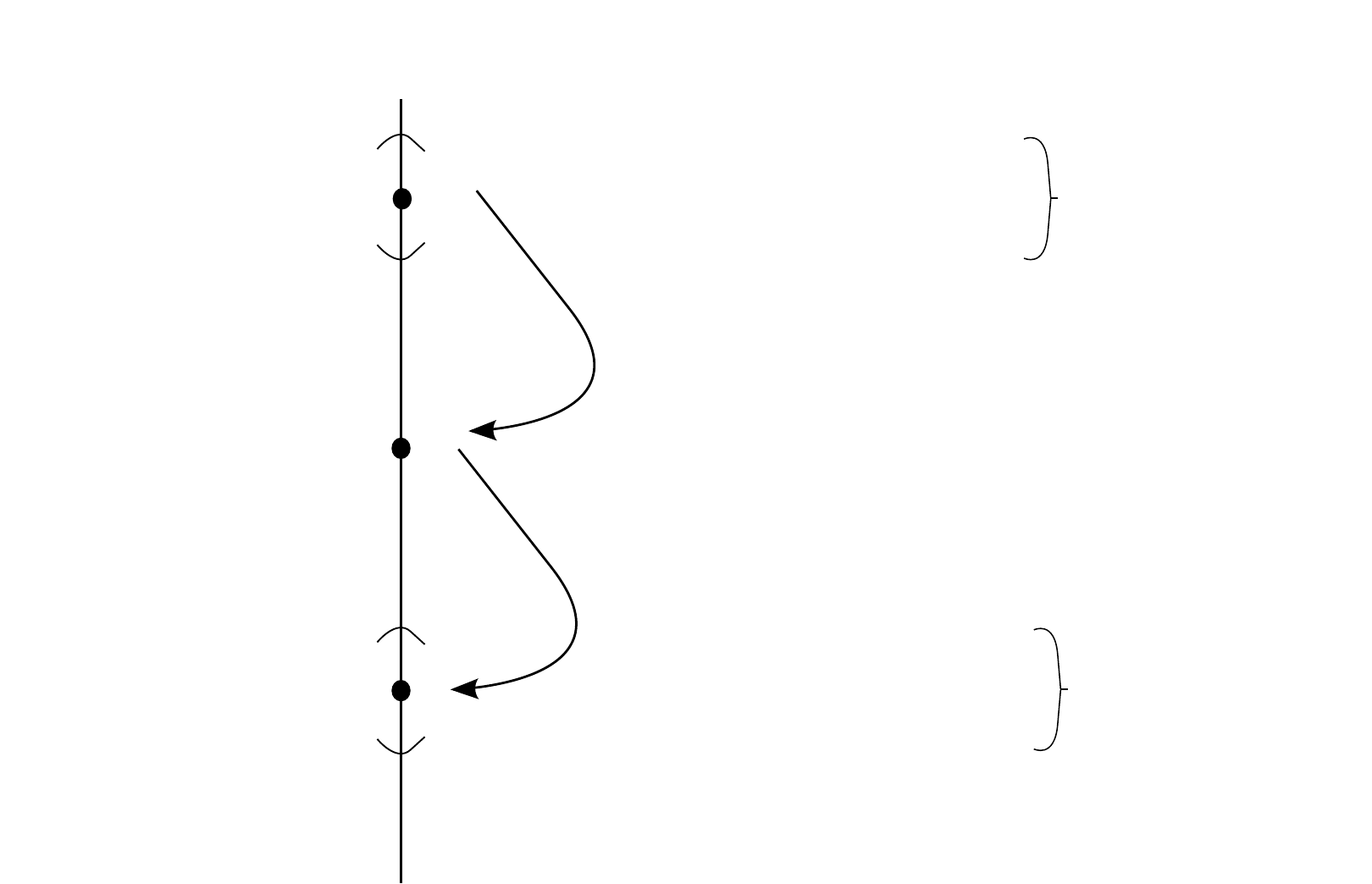
  \caption{}\label{fig:existence_y}
\end{figure}

\proof For ordinary homology classes $\alpha,\;\beta\in H_*(M)$ with $\deg(\alpha),\;\deg(\beta)<2n$ as in the statement of Theorem~\ref{maintheo}, consider $\Phi_{\beta*\alpha}([\widehat{\gamma}^k])$ as an element of the group $HFN_*^{(a,b)}(\theta^{\prime}_t)$ with $(a,b)$ as in Lemma~\ref{lemma:qgamma}. Since $\beta*\alpha=\mathfrak{q}[M]$, then, by (\ref{eqn:multip_quantum}) and (\ref{diagram:cap}) with $(c,c')=(a,b)$, we have
$$
\Phi_{\beta}\Phi_{\alpha}([{\widehat{\gamma}}^k])=\Phi_{\beta*\alpha}([{\widehat{\gamma}}^k])=
\Phi_{\mathfrak{q}[M]}([{\widehat{\gamma}}^k])=\mathfrak{q}\Phi_{[M]}([{\widehat{\gamma}}^k])= \mathfrak{q}[{\widehat{\gamma}}^k].
$$
Take $\sigma$ and $\eta$ generic cycles representing the ordinary homology classes $\alpha$ and $\beta$, respectively. The chain $\Phi_{\eta}\Phi_{\sigma}({\widehat{\gamma}}^k)$ represents the homology class $\mathfrak{q}[\widehat{\gamma}^k]$ and hence the orbit $\mathfrak{q}{\widehat{\gamma}}^k$ enters the chain $\Phi_{\eta}\Phi_{\sigma}({\widehat{\gamma}}^k)$ (by Lemma~\ref{lemma:qgamma}). Hence, (see Figure~\ref{fig:existence_y}) there is an \emph{orbit} $\widehat{y}$ in $\Phi_{\sigma}({\widehat{\gamma}}^k)$ connected to ${\widehat{\gamma}}^k$ and $\mathfrak{q}{\widehat{\gamma}}^k$ by trajectories which are solutions of (\ref{eqn:connectorbit_eq_tilde}). By the Ball-Crossing Energy Theorem~\ref{thm:ballcrossing}, (\ref{epsdelt}) and
\begin{itemize}
\item[i)] $\mathcal{A}_{\widetilde{K}}({\widehat{\gamma}}^k)=0$
\item[ii)]$\mathcal{A}_{\widetilde{K}}(\mathfrak{q}{\widehat{\gamma}}^k)=-h_0$,
\end{itemize}
we obtain
$$
-\epsilon > \mathcal{A}_{\widetilde{K}}(\widehat{y})> -h_0 +\epsilon.
$$

\qed

\subsection{Proof of Proposition~\ref{prop:example}}\label{section:prop} Let $\psi_t$ be a symplectic isotopy of $(M,\omega)$ and $\gamma$ a contractible loop in $M$.
We prove this proposition in four steps. In the first three (1a-1c) we assume $\gamma$ is constant and in the last step, (2), we consider the general case.\\
Assume that $\gamma$ is a constant loop, i.e. $\gamma(t)\equiv p$.
\begin{enumerate}
\item[(1a)] Take a Hamiltonian $H\colon [0,1]\times M\rightarrow \reals$ such that the Hamiltonian flow $\varphi^t_H$ associated with $H$ satisfies $\varphi^t_H(p)=\psi_t(p)$. Namely, take any function $H_t$ satisfying $-dH_t(\psi_t(p))=\iota_{\frac{d}{dt}\psi_t(p)}\omega.$ The point $p$ is fixed by the composition $(\varphi^t_{H})^{-1}\circ \psi_t$ and notice that the new isotopy is obtained by a Hamiltonian deformation of $\psi_t$ since $(\varphi^t_{H})^{-1}$ is the flow of some Hamiltonian usually denoted by $H^{\mbox{\small inv}}$. Denote this composition by $\overline{\psi}_t$.\\
\item[(1b)] There exists a Hamiltonian $H^{'}\colon [0,1]\times M\rightarrow \reals$ such that $\overline{\psi}_t=\varphi^t_{H^{'}}$ near the point $p$, since $\theta_t=dH^{'}_t$ near $p$ for some Hamiltonian $H^{'}$. The composition $(\varphi^t_{H^{'}})^{-1}\circ\overline{\psi}_t\equiv id$ near $p$ and it is obtained by a Hamiltonian deformation of $\overline{\psi_t}$ (and hence of $\psi_t$). Denote this composition by $\overline{\phi}_t$.\\
\item[(1c)] Take a Hamiltonian $K\colon [0,1]\times M \rightarrow \reals$ such that $p$ is a hyperbolic fixed point of $\varphi_{K}$. Then $\phi_t:=\varphi^t_{K}\circ \overline{\phi}_t$ is a Hamiltonian deformation of $\overline{\phi}_t$ and of $\psi_t$ such that $p$ is fixed by this isotopy and is hyperbolic. The reason is that, by (1b), $\overline{\phi}_t\equiv id$ near $p$.
\end{enumerate}
Consider now the general case where $\gamma$ is a contractible loop and denote $\gamma(0)$ by $p$.
\begin{enumerate}
\item[(2)]  Applying steps (1a) through (1c) to the point $p$, we obtain a symplectic path $\phi_t$ such that $p$ is a fixed point of $\phi_1$. There exists a loop of Hamiltonian diffeomorphisms $\eta_t$ such that $\eta_t(p)=\gamma(t)$ (see e.g. \cite[Section~5.1]{Gi:conley} for more details). Then $\gamma$ is a hyperbolic periodic orbit of the composition $\eta_t\circ\phi_t$ - a Hamiltonian deformation of $\psi_t$. \qed 
\end{enumerate}

\begin{Remark} Observe that with the previous construction we may create an isotopy with any prescribed linearization $d\varphi_K$ along a loop $\gamma$.
\qed
\end{Remark}

%\nocite{*}
\bibliographystyle{amsalpha}
\bibliography{references}

\end{document}